\documentclass[12pt,a4paper]{article}
 \usepackage{amsmath,amstext,amssymb,amscd}
  \usepackage[english]{babel}
 \usepackage{graphicx}

\newcommand{\Xcomment}[1]{}

\newtheorem{thm}{Theorem}[section]
\newtheorem{prop}[thm]{Proposition}
\newtheorem{lem}[thm]{Lemma}%[section]
\newtheorem{cor}[thm]{Corollary}%[section]
\makeatletter \@addtoreset{equation}{section} \makeatother

\usepackage{mathtools}

{\hfill$\qed$\medskip}

\def\qed{ \Box}

\begin{document}

\title{Stable sets of contracts in two-sided markets%\thanks{We would like to express our gratitude to the anonymous referees for their numerous useful comments.}
}
                 \author{Vladimir I.~Danilov
 \thanks{Central Institute of Economics and
Mathematics of the RAS, 47, Nakhimovskii Prospect, 117418 Moscow, Russia;
email: vdanilov43@mail.ru.}
  \and
Gleb A.~Koshevoy
\thanks{The Institute for Information Transmission Problems of
the RAS, 19, Bol'shoi Karetnyi per., 127051 Moscow, Russia, %and HSE University;
email: koshevoyga@gmail.com. }
  }

\date{}

\maketitle

To Michel Balinski on the occasion of 87-th birthday.

\begin{abstract}
We revisit the problem of existence of stable systems of contracts with arbitrary sets of contracts. We show that stable sets of contracts exists if choices of agents satisfy path-independence. We call such choice functions  Plott functions. Our proof is based on application of Zorn lemma to a special poset of semi-stable pairs. Moreover, we construct a dynamic process on the poset (generalizing algorithm Gale and Shapley) steady states of which are stable sets.

In Appendix we discuss Lehmann hyper-orders and establish a bijection between the set of Lehmann hyper-orders and the set of Plott  functions.\medskip

\end{abstract}

\section{Introduction}
There are many studies of stable sets of contracts in many-to-many {setup} (bilaterial markets), beginning from the fundamental paper by Gale and Shapley \cite{GS}. Let us mention %\textbf{only}
some of  the key papers: Kelso and Crawford \cite{KC}, Roth \cite{R}, Fleiner \cite{Fl},  Baiou and Balinski \cite{Bal}, Hatfield  and  Milgrom \cite{HM}. There are two sides of agents (men and women, students and colleges, doctors and hospitals, workers and firms, banks and clients), and agents of one of the sides may  sign contracts with the agents of the opposite side.  The amount of contracts which agents can sign is arbitrary. Baiou and Balinski \cite{Bal} introduced the notion of schedule matching which made it possible to consider, as a part of the contract not only the hiring of a particular worker by a particular firm, but also the number of hours of employment of the worker in the firm.
%{ Hatfield  and  Milgrom \cite{HM} started the study of matching with contracts with an unstructured set of contracts. The existence of a stable set of contracts has been proved, in this unstructured framework and for finite sets of contracts, in \cite{CY}.}

%\textbf
{Agents %should compare attractiveness sets of contracts,  their
preferences over set of contracts} are given by choice functions.
%defined on the set of contracts.
Such a viewpoint on agents preferences dates back to Roth \cite{R}. Roth also formulated a clever generalization of the maximization of usual preferences of  Gale and Shapley \cite{GS}: for existence of stable contracts, the choice functions have to be Plott functions, that are functions satisfying path-independence due to Plott \cite{P}.

{ Gale and Shapley’s concept of stability is a milestone of matching theory.}
There were two main approaches to establish  existence  of stable contracts: either, for a case of finite sets of contracts,  to use a variant of the Gale-Shapley algorithm, or, for infinite sets of contracts, to use the Tarski fixed point theorem  as Fleiner did in \cite{Fl}.

Here {we propose a new} approach to the existence problem using the Zorn lemma. We construct an auxiliary poset of semi-stable pairs, and maximal elements of this posets gives us stable sets of contracts. Because of the Zorn lemma, such maximal  elements exist. Thus we do not care about the  issue of finiteness of infiniteness.

Beside that, we construct a process defined on this poset of semi-stable pairs (transfinite in general), such that the steady states of this process give us  stable sets of contracts. This more refined approach allows us not only to prove the existence of  stable sets
of contracts, but also the structure properties of those stables sets.

The article is organized as follows. In Section 2, we introduce the general concept of a stable contract system. Although any number of agents is allowed in the model, it can be reduced to the case when there are only two agents. In Section 3 we recall some of the properties of Plott functions, that we use in our study of stability. In Section 5 we reformulate the concept of a stable system in the form of a stable pair introduced by Fleiner \cite{Fl}. Some weakening of this concept leads to semi-stable pairs. This is, in a sense, a key step, because the existence of stable pairs and stable contract systems is reduced to the existence of maximum elements in the poset of semi-stable pairs. In Section 6, we construct a process of  ``improving" semi-stable pairs (a variant of the Gale-Shapley algorithm), such that their steady states give us stable pairs and thus stable contract systems. This process has some good properties that allow us to establish the structural properties of a set of stable sets.

In the Appendix we give one more bijection between the set  of Plott functions and the set of Lehmann hyper-orders.

      \section{Stable contract systems}

There are two complementary sets $\mathcal F$ and $\mathcal W$ of agents (for example, firms and workers). For each pair of agents of opposite types $(f,w)$, $f\in \mathcal F$, $w\in \mathcal W$,
there is a set of contracts $C(f, w)$ which they can sign. Thus, $C(f)=\coprod_{w\in\mathcal W} C (f,w)$ is the set of contracts available for agent $f$.

Which subset of contracts will be signed, depends on preferences of the agents. Following accepted tradition in the literature, going back to Roth \cite{R}, preferences of agents are  specified by the corres\-pond\-ing choice functions. This means the following: suppose a set of contracts $X\subset C(f)$ is available to sign for an agent $f$, then $f$ chooses  a subset of $G_f (X) \subseteq X$ to sign. The preferences of agent $f$ are given by the choice function $G_f: 2^{C(f)}\to 2^{C(f)}$ that  sends $X$ to $G_f(X)\subset X$, $X\subset C_f$.

Similarly, the preferences of  agent $w$ are given by the choice function $F_w$ defined on the set $C(w)=\coprod_{f\in \mathcal F} C (f,w)$.

Denote by $C=\coprod_f C(f)=\coprod _w C (w)$ the set of all possible contracts. For the  given preferences of agents $G_f$, $f\in \mathcal F$, and $F_w$, $w\in \mathcal W$,  which subset $S\subseteq C$ could be implemented? Gale and Shapley \cite{GS} proposed the stability  concept for contract systems. Namely, two requirements have to be satisfied for a stable set $S$ of contracts.

Firstly, no one of agents wish to abandon any of the contracts from the proposed set $S$. This means that, for any agent $f$, the equality holds $G_f(S\cap C(f))=S\cap C (f)$, and, for any agent $w$,  it holds $F_w(S\cap C(w))=S\cap C (w)$.

Secondly,  there is  no pair  $(f,w)$ which wants to conclude a new contract. That is  if $c$ belongs to $C (f,w)$ and does not belong to $S$, then either $c$ does not belong to $G_f (S\cup c)$, or $c$ does not belongs to $F_w(S\cup c)$.

If both of these conditions are met, the contract system $S$ is called a \emph{stable set}.

It is convenient to rewrite these conditions in more compact form by aggregating all workers into one Worker, and all firms into one Firm. The Worker's preferences are represented by a choice function $F$ defined on the set $C$ being the union of individual ones,
                                               $$
                              F(X)=\coprod _w F_w(X\cap C(w)).
                                                 $$
Similarly, the preferences of the Firm are given by the choice function $G$ on $C$,
                                                   $$
                              G(X)=\coprod _f G_f(X\cap C(f)).
                                                     $$
Then, the subset $S \subseteq C$ is stable if and only if
\begin{itemize}
          \item[S1.] $F(S)=S$, $G(S)=S$;
\item[S2.] If $c\notin S$, then either $c\notin F (S\cup c)$, or $c\notin G (S\cup c)$.
\end{itemize}
\medskip

      Here is a couple of examples.\medskip

\textbf{Example 1.} There are two agents. The set $C$ consists of six contracts (depicted below as circles). The choice of the Worker is defined by maximization of $u_1$ {(if the utility is non-negative)} and the choice of the Firm is defined by {similar} maximization of $u_2$.

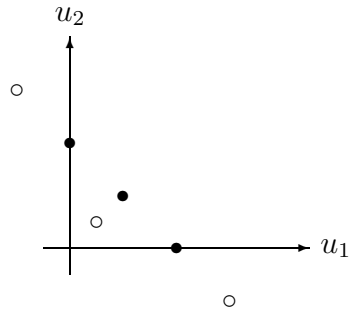
\unitlength=.7mm
\special{em:linewidth 0.4pt}
\linethickness{0.4pt}
\begin{picture}(110,64.00)(-35,0)
\put(45.00,10.00){\vector(1,0){50}}
\put(50.00,5.00){\vector(0,1){45}}
\put(50.00,30.00){\circle*{2.00}}
\put(60.00,20.00){\circle*{2.00}}
\put(70.00,10.00){\circle*{2.00}}
\put(40.00,40.00){\circle{2.00}}
\put(80.00,0.00){\circle{2.00}}
\put(55.00,15.00){\circle{2.00}}
\put(100,10.00){\makebox(0,0)[cc]{$u_1$}}
\put(50.00,54.00){\makebox(0,0)[cc]{$u_2$}}
\end{picture}

Stable contracts are depicted by black circles.\medskip

\textbf{Example 2.} The set $C$ consists of two contracts, $a$ and $b$. The preferences of the Firm is the  choice function $G$ and for the Worker it is $F$ that are specified  as follows:\medskip

      $F(a)=a$, $F(b)=\emptyset $, $F(a,b)=\{a,b\}$;

      $G(a)=a$, $G(b)=b$, $G(a,b)=b$.\medskip

We claim that there are no stable sets in such a case. In fact, the empty set $\emptyset $ is unstable, since $a$ is better for both. The singleton $a$ is unstable due to the presence of a contract $b$: $F(a,b)$ contains $b$, as well as $G (a,b)$. The set $C=\{a,b\}$ is unstable, because the Firm prefers to abandon $a$. The singleton $b$ is unstable, since the Worker refuses such a contract.\medskip

Example 2 shows that for the existence of a stable set, the choice functions $F$ and $G$ have to satisfy certain assumptions. %Let $G$ be the CF of  the Worker on the set $C$.

Choice functions, abstracted from choices of 'best variants', are of special interest in choice theory. There are two interconnected approaches to formalize the notion of a choice of 'best variants'. Due to the first one, we have to pick a preference structure on $C$, and form, with help of it, a choice operator. Due to the second (or dual) one, we impose some requirements like 'rationality' or 'consistency' on choice operators and try to build a corresponding choice mechanisms (see, for example, Aizerman and Malishevsky \cite{AM}).

Such requirements of  `rationality'  are as follows.

\begin{itemize}
\item Suppose for  a set of available contracts $X\subseteq C$ to the Firm, some contract $c\in X$ is not chosen from $X$, $c\notin G(X)$. Then for the set of available contracts $X-c$, the choice remains the same as for $X$, $G (X-c)=G (X)$. This  property of the choice functions is called \emph{independence from rejected alternatives}, or the \emph{Outcast} property.

\item Now, let,  for  a set of available contracts $X\subseteq C$ to the Firm, some contract $c\in X$  is chosen from $X$, $c\in G (X)$. Suppose the set of available contracts diminished,  for some reason, to a subset $Y\subseteq X$, but $c$ remains available in $Y$, $c\in Y$. Then $c$ has to be chosen from $Y$, $c\in G(Y)$. %\textbf
{This  property is called
the \emph{Heredity.}}
\end{itemize}

The second requirement means that contracts are \emph{substitutes} for each other, they are interesting in themselves, and not because of their connection with another contracts.

Note that the second requirement  of substitutability was violated in Example 2.

Roth \cite{R} have shown that these two properties are %\textbf
{sufficient}  for the existence of stable sets.
Aizerman and Malishevsky \cite{AM} have shown that these two requirements on choice functions are equivalent to the so-called {\em path independence}, introduced for the first time by Plott \cite{P}. Path independence means that the result of choice does not change when we divide a given set in two groups $X$ and $Y$ and then take the choice from a smaller set consisting of the group $X$ and the 'winners' of the group $Y$. Choice functions satisfying path independence are called Plott functions.

In the following section we collected results on  Plott functions that will be needed further.

           \section{Plott functions}

The Plott condition is expressed by the following compact formula:
$$
      G(X\cup Y)=G(G(X)\cup Y) .\ \ \ \ \ \ \ \ \ \ \ \ \ (PI)
$$
An equivalent condition is $G(X\cup Y)=G(G(X)\cup G(Y))$. %\textbf
{Obviously, this implies the idempotency of $G$: $G(G(X))=G(X)$ for any $X\subset C$.}

Here are some important properties of Plott functions which are useful for our analysis of stability. \medskip

\subsection{Noetherian linear orders}\label{zero}

The set of Plott functions is a huge set. They form a semi-lattice wrt union with the top element
$\mathbf 1 (X)=X$, $X\subset C$  and the minimal element is $\mathbf 0(X)=\emptyset$, $X\subset C$.

Namely, for a  family $(G_i, i\in I)$ of Plott functions, the union $\cup_i G_i$  is a Plott function as well. Let us prove this claim for two functions. If $F$ and $G$ are Plott functions on $C$, then  $F\cup G$ is also a Plott function,
which means  that the following equality holds
      $$
      (F\cup G)(X\cup Y)=(F\cup G)((F\cup G)(X)\cup Y).
      $$
The RHS is the union of $F(FX\cup GX\cup Y)$ and $G(FX\cup GX\cup Y)$. Because $GX\subseteq X$, we have  $F (FX\cup GX\cup Y)=F(X\cup GX\cup Y)=F(X\cup Y)$ . Similarly, $G(FX\cup GX\cup Y)=G(X\cup Y)$. Hence by the definition $(F\cup G)(X\cup Y)=F(X\cup Y)\cup  G(X\cup Y)$ and we get the claim proved.

This property of stability holds for any family of Plott functions. In particular, we have that aggregating Plott functions (see Section 2), yields a Plott function.\medskip

Single-valued Plott functions are in bijection  with Noetherian linear orders.\medskip

\textbf{Definition.} A linear order $<$ on $C$ is called \emph{Noetherian} if there are no infinite strictly increasing sequences $c_1<c_2<...$ of elements of $C$. Equivalently, we can say that any non-empty subset of $X\subseteq C$ has a maximum element.\medskip

For a Noetherian order $<$, let us define a choice function $G=\max_<$ by sending a set $X$ to   $G(X)$ constituted from the maximum (relative to $<$) element in $X$ (of course, $G(\emptyset )=\emptyset $). It is easy to see that such a choice function $G= \max_<$ satisfies the Plott condition. Note that without the Noetherian condition, the statement is incorrect.

Let us note one characteristic property of the choice function $G= \max _<$.
Such a choice function is ``single-valued": $G (X)$ contains exactly one element (if, of course, $X\ne \emptyset $). The converse is also true: if $G$ is a single-valued Plott function, it is given by maximization of some Noetherian linear order $<$. This order $<$ can be given explicitly: for different $a, b\in C$
                $$
      a<b \text{ iff } G(a,b)=b.
                  $$
It is easy to understand that the relation $<$ is transitive and complete, so that we get a linear order. As for the Noetherian property:  suppose there is an infinite increasing sequence $a_1<a_2<...$. Let us form the set $A=\{a_1,a_2,...\}$. $G (A)$ consists of a single element; let this element be $a_n$. Due to  the heredity property (see \ref{two}) for  $\{a_n, a_{n+1}\} \subseteq A$, we get $G(a_n,a_{n+1}) =G(A)=a_n$. On the other hand, from $a_n<a_{n+1}$ we see that $G(a_n, a_{n+1})=a_{n+1}$. A contradiction.\medskip

According to the Zermelo theorem, there is a plethora of Noetherian linear orders on $C$, so we get  a plethora of Plott functions.

One can construct a non-empty valued Plott functions from Noetherian linear orders.
For a Plott function $G$, let us take the set of Noetherian linear orders $<_i$ , $i\in I$, which are inferior to $G$. That is, for any $X\subset C$, it holds $\max_{<_i}(X)\subset G(X)$, $i\in I$.
Then $G$ is equal to the union  of Plott functions $\max_{<_i}$, $i\in I$. This construction of Plott functions generalizes the theorem of Aizerman and Malishevski, and was established in \cite{DKS}.

\subsection{Heredity and Outcast}\label{two}  From the definition of a Plott function, the following inclusion holds
      $$
          G(X\cup Y)\subseteq G(X)\cup G(Y).    \ \ \ \ \ \ \ \ \             (S)
      $$
For sets  $A\subseteq B$, let us set $X=A$ and $Y=B\setminus A$. Then from (S) we get
\[
G(B)\subset G(A)\cup G(B\setminus A),
\]
that is nothing but the {\em heredity } requirement (\cite{AM}): for $A\subseteq B$
      $$
                       G(B)\cap A\subseteq G(A).
                                 $$
%\textbf
{In other words, if an element $a$ belongs to a small set $A$ and is chosen from a bigger set $B$ then it has be chosen from the smaller set $A$. }

Note that inclusion (S) is valid for any number of sets: \emph{for any family $(X_i, i\in I)$ of subsets of $C$, there is an inclusion}:
$$
        G(\cup _i X_i)\subseteq \cup _i G(X_i).
      $$

%\subsection{Outcast}\label{three}.

%\textbf
{The second important property of Plott functions is} the \emph{outcast} property: assume that $G(X)\subseteq Y\subseteq X$; then $G (Y)=G (X)$.\medskip

Indeed, $G (Y)=G(G(X)\cup Y)=G(G(X)\cup G(Y))=G(X\cup Y)=G(X)$.\medskip

%\textbf
{%Conversely,
These two properties characterise Plott functions: if a choice function $G$ satisfy Heredity and Outcast, then $G$ is a Plott function.}

\subsection{Closure operator}\label{five}

Note that, for a Plott function, the property (PI)  is valid not only for two sets, but for any number of them. More precisely, for any family  $(X_i, i\in I)$ of subsets of $C$, the following equality holds
      $$
                            G(\cup _i X_i)=G(\cup _i G(X_i)).
      $$
Indeed, if $Y=   \cup _i G(X_i) $ and $X= \cup _i X_i$, we have (see \ref{two})
$$
G(\cup _i X_i)\subseteq  \cup _i G(X_i) \subseteq \cup _i X_i.
      $$
%\textbf
{Due to Outcast}, we have  $G(\cup _i X_i)=G(\cup _i G(X_i))$.

%\textbf
{ By applying the above equality %this property
to the collection $\mathcal I(X)=\{Y\subseteq C, \ G(Y)=G(X)\}$, we get that this collection has the largest element, namely the union $Z=\cup_{Y\in \mathcal I(X)} Y$. Indeed,
}$$
G(Z)=G(\cup_{Y\in \mathcal I(X)} Y)=G(\cup_{Y\in \mathcal I(X)} G(Y))=G(G(X))=G(X).
$$
%\textbf
{We denote by $G^*(X)$ this largest set.
In other words, for any $X\subseteq C$, the following two statements are equivalent:}

(1) $G(Y)=G(X)$,

(2) $X\subseteq Y\subseteq G^*(X)$.

Thus, for a Plott function $G$, we get an expanding operator $G^*$. It can be shown that $G^*$ is a closure operator.  Moreover, for a Plott function $G$,
the closure $G^*(X)$ is a kind of convex hull of $X$ (for details see \cite{K} for finite $C$ and  \cite{MP} in the general case). The inversion to the operator $G^*$ is defined by the formula:
                     $$
                              G(X)=\{x\in X, x\notin      G^*(X-x)\}.
                       $$

In particular, when $X=\emptyset $, we get that there is the largest subset of $G^*(\emptyset )$, the choice from which is empty. It can be called \emph{Nil-set} for a Plott function $G$, $Nil(G)$ (or   Dummy set, the complement to it is the `support' of $G$). This nil-set does not affect the choice in any way: for any $X$, $G(X)=G(X-Nil (G))=G(X\cup Nil (G))$ holds. \medskip

\subsection{Blair hyper-order}

To compare subsets in $C$ in terms of their attractiveness for a G-agent with Plott function $G$, Blair \cite{B} proposed %\textbf
{to use} a natural (hyper)-relation $\preceq _G$ on $C$. Namely, for subsets $A$ and $B$ of $C$, we set
                                 $$
       A\preceq _G B, \text{ if } G(A\cup B)\subseteq B.
      $$
In this case, $G (A\cup B)=G(G(A\cup B)\cup B)=G(B)$. In other words, $A\preceq_G B$ is equivalent to the inclusion $A\subseteq G^*(B)$ (as well as to the inclusion $G^*(A)\subseteq G^*(B)$). %Meaningfully,
The relation $A\preceq _G B$ means that adding $A$ to $B$ does not affect the preferred elements in $B$, and therefore, adding $A$ to $B$ does not increase the attractiveness of $B$. %\textbf
{For any $A$, it holds $G^*(A)\preceq _G G(A)$.}

%For the $F$-agent, the hyper-relation $\preceq_F$ is defined similarly: $A\preceq_F B$ if $F(A\cup B)\subset B$.

\begin{lem}\label{lemma1}
{The hyper-relation $\preceq =\preceq _G$ is transitive.}\end{lem}
{\em Proof}.
Let  $A\preceq_G B$ and $B\preceq_G D$. Then %$A\preceq D$. In fact,
$G (A\cup D)=G (A\cup G(D))=G(A\cup G(B\cup D))=G(A\cup B\cup D)=G(G(A\cup B)\cup D)=G(G(B)\cup D)=G(B\cup D)=G(D)$. That is nothing but $A
\preceq_G D$. \hfill $\Box$\medskip

Another important property of this hyper-relation is its consistency with the unions.

%\textbf{Lemma 2.} \emph
\begin{lem}\label{lemma2}
{Let $(X_i, i\in I)$ be  a family of  subsets in $C$, and $X_i\preceq _G Y$ for any $i\in I$. Then $\cup _i X_i\preceq _G Y$.}\end{lem}

The consistency shows the ``qualitative" character of the preorder $ \preceq_G$: the quantity can not bit
the quality.
%attractiveness of the union is not higher than the attractiveness of the terms.

{\em Proof}. Let $X=\cup _i X_i$. Then
                                       $$
      G(X\cup Y)=G((\cup _i X_i)\cup Y)=G(\cup _i (X_i\cup Y))=G(\cup _i
      G(X_i\cup Y))=G(\cup _i G(Y))=G(Y).
                                         $$ %\hfill $\Box$ \medskip

%Another important property (the `frame') is given in \cite{DKS}, but we won't need it.\medskip

      \section{Properties of stable contract systems}

From now on we suppose that \emph{the choice function $G$ of the Worker  and the choice function $F$ of the Firm are two Plott functions} defined on the set $C$ of contracts.
For given Plott functions $F$ and $G$, let us denote the Worker as $F$-agent and the Firm as $G$-agent.

Recall that a \emph{stable contract system} (or a \emph{stable set}) is a subset $S \subseteq C$ such that S1 and S2 are satisfied. Note that  S2 can be rewritten as follows
\begin{itemize}
          \item[S2'.] $F^*(S)\cup G^*(S)=C$.
          \end{itemize}

%\textbf
{To see that, let $c\in F (S\cup c)$; then $c\notin F^*(S)$. Symmetrically, $c\notin G^*(S)$. Which violates S2'. Conversely, let $c$ not belong to $F^*(S)$. Then $c\in F(S\cup c)$, because otherwise from Outcast we would have $F (S\cup c)=S$ and $c\in F^*(S)$. Symmetrically, if $c$ does not belong to $G^*(S)$, then $c\in G (S\cup c)$, contrary to S2.}\medskip

Denote by {\bf ST} the set of stable subsets in $C$. We are interested in the structure of this set. %{\bf ST}.
Firstly, we prove that {\bf ST} is nonempty and, as a rule, contains several stable sets. Secondly, %\textbf
{we show that if one stable system is better than another for the $F$-agent then it is %worse
vice versa for the $G$-agent. Namely, the following assertion is true.}

%\textbf{Proposition 1.} \emph
\begin{prop}\label{pro1}
{The restriction of %\textbf
{the Blair hyper-preorder} $\preceq _G$ to the set {\bf ST} (as well as the restriction of $\preceq _F$) is an order (that is, an antisymmetric relation).}\end{prop}

{\em Proof}.  Let $S$ and $T$ be stable sets and  $S\preceq _G T$. Then $G (S\cup T)=G (T)=T$.
If  $T\preceq _G S$ holds as well, then $G(S\cup T)=G(S)=S$, that is ,the equality $S=T$ holds. \hfill $\Box$ \medskip

A remarkable fact (\cite{B, R}) is that \emph{the orders $\preceq_F$ and $\preceq _G$ on {\bf ST} are opposite to each other.}  %In other words, if one stable set $S$ is better than $T$ for $G$-agent, then $T$ is better than $S$ for $F$-agent.
That is
\[
T\preceq_G S \mbox{ if and only if } S\preceq_F T.\]

Roth \cite{R} called this as \emph{the polarization of the interests of the opposite parties}. The proof is based on the following lemma (see \cite[Lemma 20]{Larxiv}).
%\textbf{Lemme 3} (see \cite[Lemma 20]{Larxiv}). \emph
\begin{lem}\label{lemma3}
{Let $S$ be a stable set, and $T$ be an arbitrary subset in $C$. If $S\preceq _G T$, then $G(T)\preceq_F S$.}\end{lem}

We give a more concise proof than \cite[Lemma 20]{Larxiv}.

{\em Proof}. If $S \preceq_G T$, then $G(S\cup T)=G(T)$. Let $x$ be an arbitrary element of $G(T)-S$. We state that $x\preceq_F S$. From this, obviously, follows $G(T)\preceq_F S$.

Thus, we have to check the relation  $x\preceq _F S$, that is  $x\notin F(S\cup x)$. On the contrary, suppose $x\in F(S\cup x)$. Then, due to S2, we have $x\notin G(S\cup x)$. Due the heredity of $G$, we get $x\notin G(S\cup T)=G(T)$, which contradicts the fact that $x\in G(T)$. \hfill $\Box$\medskip

 %\textbf
 {Furthermore,} we consider the set {\bf ST} as a poset with the partial order $\preceq =\preceq _G$.
 We will see below that this poset is a nonempty complete lattice.
 Knuth is credited for the discovery of this remarkable fact which %refers to was voiced by Knut and was
established in full generality by Blair \cite{B}. In particular, there is a stable set that is least attractive for the $F$-agent
and most attractive for the  $G$-agent.
We postpone the construction to Section 6, where we get this fact as the result of some natural process generalizing the Gale-Shapley algorithm.

\section{Stable and semi-stable pairs}

Since $F$ and $G$ are Plott functions,
 we can reformulate (see \cite{Fl}) the stability condition.\medskip

\noindent \textbf{Definition.} \emph{A stable pair} is a pair $(Y,Z)$ of subsets in $C$ such that the  following properties are
satisfied:
\begin{itemize}
\item[SP1] $Y\cup Z=C$,
\item[SP2] $G(Y)=F(Z)$.
\end{itemize}

This means a division of competencies. Namely, the entire set of contracts $C$ is covered by  two sets $Y$ and $Z$, the $G$-agent chooses from $Y$, and the $F$-agent chooses from $Z$. If their choices coincide, we can expect that they %give
yield a stable set. The following Lemma states the equivalence of stable pairs and  stable sets.\medskip

%\textbf{Lemma 4.} \emph
\begin{lem}\label{lemma4}{If $(Y,Z)$ is a stable pair, then $S=G (Y)=F (Z)$ is a stable set. Conversely, if $S\subseteq C$ is a stable set, then the pair $(G^*(S), F^*(S))$ is a stable pair.}\end{lem}

{\em Proof}. Let  $S$ be a stable set, $Y=G^*(S)$ and $Z=F^*(S)$. Then,  since $F(Z)=F (F^*(S))=F (S)=S$ and $G (Z)=G (G^*(S))=G(S)=S$, we get the property SP2. The SP1 is nothing but S2'.  %property was checked above.

Conversely, let the pair $(Y,Z)$ be stable and $S=G (Y)=F (Z)$. The property S1 follows from the idempotence of $F$; indeed, $F (S)=F(F(Z))=F(Z)=S$. Similarly, $S=G (S)$. The property S2' follows from the fact that $Y\subseteq G^*(S)$ and $Z\subseteq F^*(S)$. $\Box$ \medskip

If the pair $(Y,Z)$ is stable, we call $S=G (Y)=F (Z)$ a stable set \emph{corresponding to the pair} $(Y, Z)$. But how to construct stable pairs? To do this, we will slightly weaken the concept of stability by introducing a degree of asymmetry of agents.\medskip

\noindent \textbf{Definition.} A pair $(Y,Z)$ of subsets in $C$ is called \emph{semi-stable} if the following conditions hold

\begin{itemize}
\item[{\bf SSP}1] $Y\cup Z=C$,
\item[{\bf SSP}2] $G(Y)\subseteq F(Z)$.
\end{itemize}

Of course, if {\bf SSP}2 is satisfied as an equality, then we get validity of SP2 and, hence,  $(Y,Z)$ is a stable pair. On the other hand, it is easier to find pairs which satisfy the semi-stability conditions; for example, the pair $(\emptyset, C)$ is semi-stable. Denote the {\bf SSP} set of semi-stable pairs, and endow it with the following order relation:
$$
(Y,Z)\le (Y', Z'), \text{ if } Y\subseteq Y', Z'\subseteq Z.
        $$
In other words, the first component increases, and the second decreases.

For a semi-stable pair $(Y,Z)$, define a new pair $\Phi (Y,Z)=(Y',Z')$, where
\begin{equation}\label{psi}
      Y'=Y\cup F(Z), \quad Z'=Z-F(Z)\cup G(F(Z)).
\end{equation}

\begin{lem}\label{A}
For  a semi-stable pair $(Y,Z)$, the pair $(Y', Z')$ defined in (\ref{psi}) is also  semi-stable.
\end{lem}

Indeed,  {\bf SSP}1 obviously holds,  since we have   $Y'\cup Z'=Y\cup Z\cup G(F(Z))=C$.
It remains to check that $G(Y')$ is contained in $F(Z')$. Because $G(Y)\subseteq F (Z)$, we have  $G(Y')=G(Y\cup F(Z)) =G(G(Y)\cup F(Z))=G(F(Z))$. Therefore, we have to check that $G(F(Z))$ is a subset of $F(Z')$. Due to the heredity of  $F$, from the inclusion of $Z'\subseteq Z$, we obtain the inclusion of $F(Z)\cap Z'\subseteq F(Z')$. Since $G(F(Z))$ is a subset of  $Z'$ and subset of $F(Z)$, we conclude $G(F(Z))\subset F(Z)\cap Z'\subseteq F (Z')$. \hfill $\Box$

%\textbf{Proposition 2.} \emph
\begin{cor}\label{pro2}{Maximal elements of the poset $({\bf SSP},\le )$ are stable pairs.}
\end{cor}
{\em Proof}. Let a pair $(Y,Z)$ be a maximal element of the poset $({\bf SSP},\le )$. Since $(Y, Z)\le (Y', Z')$, it holds
that  $Z=Z'$, and hence, by the definition $Z'$, we have $F(Z)=G(F(Z))$. Because $Y=Y'$, $F(Z)\subset Y$, and we get the inclusions
$G(Y)\subset F(Z)\subset Y$. Therefore, from the outcast property for $G$, we get $G(Y)=G(F(Z))$.
Thus $G (Y)=F (Z)$ and the pair $(Y,Z)$ is stable.  \hfill $\Box$\medskip

We claim that \emph{the poset $({\bf SSP},\le) $ is an inductive ordered set} (or {\em chain complete poset}), which means that any chain in this set has an upper bound. Then  maximal elements in the poset $({\bf SSP},\le) $ exists by Zorn's Lemma. Moreover, for any pair $(Y,Z)$ there is a maximal pair that is superior to $(Y,Z)$. %\textbf
{
In order to establish %inductivity
chain completeness of the poset $({\bf SSP},\le) $ we prove a more general statement.}

%\textbf{Lemma 5.} \emph
\begin{lem}\label{B}{Let $((Y_i, Z_i), i\in I)$ be a family of semi-stable pairs, $Y=\cup _i Y_i$, $Z=\cap _i Z_i$. Then the pair $(Y,Z)$ is semi-stable.}\end{lem}

{\em Proof}. Because of the inclusion $Y\cup Z_i\supseteq Y_i\cup Z_i=C$ and
the distributive law, we have
$$
      Y\cup Z=Y\cup  (\cap _i Z_i)=\cap _i (Y\cup  Z_i)=\cap _i C=C.
                                        $$
%\textbf
{To verify} the inclusion $G (Y)\subset F (Z)$, firstly, let us establish the inclusion $G(Y)\subseteq Z$. Suppose  $y\in G (Y)$ and $y\not\in Z$. In such a  case, $y$ does not belong to some $Z_i$.  From $Y_i\cup Z_i=C$ follows that $y\in Y_i$. From the heredity of $G$, we get that $y\in G (Y_i)$, hence $y\in F(Z_i)\subseteq Z_i$. This contradiction
proves the inclusion $G(Y)\subseteq Z$.

Let $y$ be an element of $G(Y)$ and suppose that for some $i$,  $y\in G (Y_i)$.
Then, because $G(Y_i)\subseteq F(Z_i)$,  we get $y\in F(Z_i)$.
Since $y\in Z$,  from  the heredity of $F$  we get  $y\in F (Z)$. \hfill $\Box$\medskip

Thus the poset $({\bf SSP}, \le) $ has a maximal element. Hence,  due to Corollary \ref{pro2}, such an element is a stable pair and gives a stable set.

\section{Process of sequential improvement}

In the previous section, we have shown that, for any semi-stable pair $(Y,Z)$, there exists a stable set $S$ which is an upper bound to $(Y,Z)$. In fact,  we can get more than just an existence theorem. Namely, we can us the idea of the above proof to define a process of construction of a stable set that is not only an upper bound to the original pair, but also is the minimal upper bound
(in the sense that will be defined later).

Specifically, for each semi-stable pair $(Y,Z)$, we defined by (\ref{psi}) a new semi-stable pair $\Phi (Y,Z)=(Y',Z')$. Thus, on the set {\bf SSP} of semi-stable pairs, the monotonic dynamics $\Phi :{\mathbf{SSP}}\to{\mathbf {SSP}}$ is defined.
The fixed points of this dynamics are stable sets. The transformation $(Y,Z)$ to $\Phi (Y,Z):=(Y',Z')$ is nothing else but  the Gale-Shapley algorithm: at the step $(Y,Z)$, %the Worker,
the $F$-agent, makes  an offer $F(Z)$ to % the Firm,
the $G$-agent, and than the $G$-agent accepts the part $G(F(Z))$ of the offer and rejects the rest $F(Z)-G(F(Z))$. %After that,
The process continues at the step $(Y',Z')$ and so on.\medskip

\textbf{Remark.} %The main role in forming $(Y', Z')$ from the pair $(Y,Z)$ plays $Z$. If $Z$ stabilizes, we get a stable pair.
Lehmann (\cite{Larxiv}) %used exactly it, writing out
defined the dynamics only in terms of $Z$, $Z'=Z-(F(Z)-G(F(Z))$. Due to this %However, with this
approach, it is not very clear which $Z$  we  can take as the initial state.
Lehmann showed that if we begin the process with $Z=C$, it stabilizes and gives a stable set. Note that our process can be started with any semi-stable pair. This flexibility to choose the initial pair of the process allows us to get       some interesting lattice properties of the set of stable pairs.\medskip

Let $P_0= (Y_0, Z_0)$ be a starting  semi-stable pair. We get pairs $P_1=\Phi (P_0)$, $P_2=\Phi (P_1)$ and so on. However, this process may not lead to a stable pair in a finite number of steps, it may run indefinitely. In such a case, let us  define $P_\omega $ to be the pair $(\cup _k Y_k, \cap _k Z_k)$ and  restart  the process at the pair $P_\omega $. That is we have to consider a transfinite  process.

This means that, for any ordinal number $ \alpha $ (that is, for any %\textbf
{well}-ordered set $\alpha $), we must define a semi-stable pair $P_\alpha =(Y_\alpha, Z_\alpha )$. This is done in two different ways, depending on whether the number $\alpha $ is a limit or not.
A number $\alpha $ is called \emph{limit} if it has no immediate predecessor, that is, it does not have the form $\alpha =\beta +1$. For a non-limit number  $\alpha =\beta +1$, we set $P_\alpha =\Phi (P_\beta)$, and,  for a limit number  $\alpha$, we set $P_\alpha =(\cup _{\beta <\alpha }Y_\beta, \cap _{\beta <\alpha }Z_\beta)$.

Because of lemma \ref{B}, for any ordinal $\alpha$, the pair $P_\alpha $ is semi-stable. Since the dynamics of $\Phi $ is monotone, sooner or later the transfinite sequence of $P_\alpha $ reaches the steady state. % it is stabilizing.
Let us check that the decreasing sequence $Z_\alpha $ reaches the steady state %stabilizes
(then the sequence $Y_\alpha$ also stabilises).
Indeed, consider the increasing sequence $\hat Z_\alpha :=C-Z_\alpha $. If,  for each step of  transition from $ \alpha $ to $\alpha +1$, the set of $\hat Z_\alpha $  strictly increases, then the cardinality of $\hat Z_\alpha $ is not less than the cardinality of $\alpha $. Hence we get a contradiction at the step when the cardinality of $ \alpha $ is greater than the cardinality of $C$.

Thus, we begin the transformation process at an arbitrary semi-stable pair $P_0$
and get a (transfinite) sequence $(P_\alpha )$, which has a steady state for large $\alpha $. This steady state is a stable pair denoted by $\Phi _\infty (P_0)=(Y_\infty, Z_\infty)$. According to the proof of Corollary \ref{pro2}, this pair is stable and gives the corresponding stable system of contracts, which we denote by $\sigma (Y_0, Z_0)=S$. Since the process $P_\alpha $ improves the position of the $G$-agent all the time (and the position of the $F$-agent is getting worse), we get that $Y_0 \preceq _G S$ (respectively, $S\preceq _F Z_0$). So the result of the process is a stable set $S$, which is no worse for the $G$-agent than the initial state $Y_0$.

In fact, the resulting stable set $S$ is minimal with the property of being better than $Y_0$. This follows from

%\textbf{Proposition 3.} \emph
\begin{thm}\label{pro3}
{Let $(Y, Z)$ be a semi-stable pair, $T$ be a stable set, and $Y\preceq _G T$. Then, for the limit stable set $S=\sigma (Y,Z)$,
it is true that  $S\preceq _G T$.} \end{thm}

{\em Proof}. $Y\preceq _G T$ is equivalent to the inclusion $Y\subseteq G^*T$. We claim that $T \subseteq Z$.
To show this, suppose that some $t\in T$, and $t$  does not belong to $Z$. Since $Y\cup Z=C$, $t\in Y$ and, hence, $t\in G^*(T)$. Because $T$ is a stable set, we have $G (G^*(T))=G(T)=T$, that implies the inclusion $t\in G (G^*(T))$. From the heredity property of $G$, due to $t\in Y\subseteq G^*T$, we have   $t\in G(Y)$. Due to S2 $G(Y)\subseteq F(Z)$. Therefore, $t\in Z$,  a contradiction.

Therefore $T\subseteq Z$, and hence $T\preceq _G F (Z)$. Because of Lemma \ref{lemma3}, we have $F(Z)\preceq _G T$. Since $Y\preceq _G T$, for $Y'=F\cup F(Z)$, we get (see the proof of Proposition \ref{pro2}) that $Y'\preceq _G T$. That is, for the pair $(Y', Z')=\Phi(Y,Z)$, we have the same relation $Y'\preceq_G T$. Hence,  for any pair $(Y_\alpha, Z_\alpha )$, we have $Y_\alpha \preceq _G T$. Therefore, for the limit pair $(Y_\infty, Z_\infty )$, we have $Y_\infty \preceq _G T$. Since $S\subseteq Y_\infty $, we conclude $S\preceq _G T$.    \hfill  $\Box$

      \section{Some consequences}
Here are some consequences of our results.

Theorem \ref{pro3} shows that, for given initial semi-stable pair $(Y, Z)$,  the steady state $S=\sigma (Y,Z)$ of the process $\Phi $ is not  only a stable set, but a minimal stable  state among the set of stable sets which dominate $Y$. In particular, for the initial semi-stable  pair $(\emptyset, C)$, the stable set  $\sigma (\emptyset, C)$ is the worst for $G$-agent and the best for $F$-agent, a fact discovered in \cite{GS}.

Moreover, \emph{the poset $({\bf ST}, \preceq_G)$ is a complete lattice}.
In other words, for any family $(S_i, i\in I)$ of stable sets, there is the least upper bound $S=\vee _i S_i$. Namely, due to Lemma \ref{B}, for the pairs $(G^*(S_i), F^*(S_i))$, we get the semi-stability of the pair $(Y,Z)$, where $Y=\cup _i G^*(S_i)$, $Z=\cap _i F^*(S_i)$.
Let $S=\sigma (Y,Z)$ be the corresponding stable set, then from Theorem \ref{pro3},  for any $i\in I$, it follows that $S_i\preceq _G S$ holds. On the other hand, if $T$ is a stable set and, for any $i$,
$S_i\preceq _G T$, then $G^*(S_i)\preceq _G T$, and due to Lemma \ref{lemma2}, we get
$\cup _i G^*(S_i)\preceq _G T$. Therefore, due to  Theorem \ref{pro3} $S\preceq _G T$. But this means exactly  that $S$ is the least upper bound for the family $(S_i)$.

Another consequence% already on the topic of
is related to comparative statics. Suppose that $S$ is a stable set with respect to agents with choice functions  $F$ and $G$. Suppose that the preferences of the $F$-agent have changed to a new choice function $F'$, such that  $F\le F'$, that is, for any $X\subset C$, $F(X)\subset F'(X))$.  In other words, the $F$-agent has ``weakened"  requirements for the best contracts.
Note that if $S\subset C$ is a stable set for agents with choice functions $F$ and $G$, $S$ remains stable for agents with choice functions $F'$ and $G$.
%Of course, the old stable set $S$ will cease to be stable.

For a stable set $S$, consider a pair $(Y', Z')$, where $Y'=G^*(S)$ and  $Z'=F'^*(F^*(S))$.
 $(Y', Z')$ is a semi-stable pair with respect to agents with choice functions $F'$ and $G$.
It is easy to see that {\bf SSP1}, $Y'\cup Z ' =C$, holds true, and {\bf SSP2} follows as well,
                                                             $$
                       G(Y')=S=F(F^*(S))\subseteq       F'(F^*(S))=F'(Z').
                                                               $$
Let $S'=\sigma '(Y', Z')$ be the corresponding limit stable for the process $\Phi$ defined with respect to  to $F'$ and $G$.  Then, according to Theorem \ref{pro3}, $S'$  is not worse than the initial $Y'=G^*(S)$, $Y'\preceq _G S'$, that is, $S\preceq _G S'$.

%We see that in this case there is
This defines a natural transition from the old stable sets $S$ to the new $S'$, $S\preceq _G S'$. In other words, weakening the requirements of the $F$-agent improves the position of the $G$-agent. From Lemma \ref{lemma3}, we see that $S'=G(S')\preceq _F S$ holds.
That is changing from $F$ to an upper $F'$, $F\le F'$, only worsening the outcome for the $F$-agent.

\section*{Appendix. Lehmann hyper-orders}

With any Plott function $G$, the Blair hyper-order $\preceq _G$ is associated. Recall that it is defined by
              $$
      A\preceq _G B \text{ iff } G(A\cup B)\subseteq B.
                                       $$
%For a Plott function $G$,
 Lehmann \cite{L}  defined another (more strong) hyper-order $\prec _G$ by
                                  $$
A \prec _G B \text{ iff } G (B)\ne \emptyset \text{ and } G(A\cup B) \text{ does not intersect with } A.
                               $$
Note that in this case, $G (A\cup B)$ is contained in $B$ and moreover coincides with $G(B)$.
In fact, we have  $G(A\cup B)=G(G(A\cup B)\cup B)=G(B)$.

This `strict' hyper-order is interesting because it is uniquely defines the initial Plott function $G$. To do this, we need to recover the set $G(B)$ for any $B\subseteq C$ in terms of $\prec =\prec _G$. Here it is convenient to introduce two notions related to the hyper-relation $\prec $. Namely, a set B is called \emph{essential} if $\emptyset \prec B$, and \emph{insignificant} otherwise.

           If $B$ is insignificant, then we have $G(B)=\emptyset $.

If $B$ is essential, then there is $G(B)\ne \emptyset $, then $G (B)=B-A$, where $A=\{a\in B, a\prec _G B\}$. In fact, $\{a\}$ does not intersect with $G(a\cup B)=G(B)$, so $a$ is not contained in $G(B)$. On the other hand, if $b\in G (B)$, then the relation $b\prec _G B$ does not performed.

Below we will get the necessary and sufficient conditions on a hyper-order $\prec$, which has the form $\prec _G$ for some Plott function $G$.

      The following are the important properties of the hyper-relation $\prec_G$.\medskip

      \textbf{Proposition A1.} \emph{The hyper-relation $\prec =\prec _G$ has following properties:}

           L0. $\prec $ \emph{irreflexive.}

           L1. Left weakening. \emph{If $A'\subseteq A\prec B$, then $A' \prec B$.}

           L2. Union. \emph{Let $(A_i, i \in I)$ be a nonempty family of subsets in $C$. If $A_i\prec B$ for any $i\in I$, then $\cup _i A_i \prec B$.}

           L3. Right strengthening. \emph{If $A\prec B\subseteq B'$, then $A\prec B'$.}

           L4. Cancelation. \emph{If $A\prec A\cup B$, then $A\prec B$.}

           L5. Domination. \emph{Any essential set dominates (in the sense of $\prec$) any non-essential one.}\medskip

\emph{Proof.} L0 is obvious.

Let's check L1. We are given that $G (A\cup B)\subseteq B$ and even $G(A\cup B)=G (B)$. Then $G (A'\cup B)=G(A'\cup G(B))=G(A'\cup G(A\cup B))=G(A'\cup A\cup B)=G(A\cup B)$ is contained in $B$.

Check L2. $G ((\cup _i A_i) \cup B)=G(\cup _i (A_i\cup B))=G(\cup _i G(A_i\cup B))$ is contained in $\cup _i G (A_i\cup B)$ and especially in $B$.

Let's check L3. $G(A\cup B')=G(A\cup B\cup B')=G(G(A\cup B)\cup B')=G(G(B)\cup B')$ is contained in $B'$, since $G(B)\subseteq      B\subseteq B'$.

Let's check L4. $G(A\cup B)=G((A\cup B)\cup B)=G(G(A\cup B)\cup B)=G(G(B)\cup B)=G (B)$ and is nonempty.

Let's check L5. Let $A$ be essential, and $B$ be insignificant. The essentiality of $A$ means that $G(A)\ne \emptyset $. The insignificance of $B$ means that $G(B)=\emptyset $. We have to check that $G (A\cup B)=G(A\cup G(B))=G (A)$ does not intersect with $B$. Let's assume that it intersects, and $b\in B$ belongs to $G(A)$. Then $b$ is selected in the larger set $A\cup B$ and belongs to the smaller one $B$; from the heredity we have $b\in G(B)$. This contradicts to the emptiness of $G(B)$. $\Box$\medskip

\textbf{Definition.} \emph{The Lehmann hyperorder} is a hyper-relation, which satisfies the properties L0-L5.\medskip

{\bf Remark}. Lehmann \cite{L} called a {\em  qualitative measure} a transitive hyper-relation which satisfies L0-L4 and such that the union of negligible sets is negligible. The next lemma shows that  the transitivity is a redundant requirement.\medskip

      \textbf{Lemma A1.} \emph{The Lehmann hyper-order is transitive.}\medskip

{\em Proof}. Let $A\prec B$ and $B\prec C$. Then $A\prec B\cup C$ (according to L3), as well as $B\prec B\cup C$. Hence by L2, it holds true  $A\cup B\prec B\cup C$, and due to L3  we have $A\cup B\prec A\cup B\cup C$. Then from L4 follows $A\cup B\prec C$, and hence from L1 $A\prec C$. \hfill $\Box$\medskip

This lemma justifies the use of the term hyper-order.

Now we construct a mapping from the set of Lehmann hyper-orders to the set of Plott functions.

Let $\prec$ be a Lehmann hyper-order and $D$ be the set of all negligible elements in $C$. For a set $A\subseteq C$, define
                              $$
      L(A)=D\cup \{c\in C, c\prec A\}.
                                 $$
Note that, for  a negligible set $A$, we have $L (A)=D$.
Because of the union property L2,  for an essential set $A$, it holds that $L (A)\prec A$.
Note that $L (A)$ is the largest subset with the property $L(A)\prec A$.\medskip

      \textbf{Lemma A2.} \emph{If $A\subseteq L (A)\cup B$, then $L(A)\subseteq L(B)$.}\medskip

{\em Proof.} For a negligible set $A$, the statement is obviously true. Let $A$ be essential, then $L (A)\prec A$. Since $A\subseteq  L(A)\cup B$, then, from L3,  we have $L(A)\prec L(A)\cup B$. From L4, we get $L(A)\prec B$. This implies that $B$ is essential and $L(A) \subseteq L(B)$. \hfill $\Box$\medskip

      Now let's us define a choice function $T=T_\prec$ by the rule
                             $$
                                   T(A)=A\setminus L(A).
                               $$

      \textbf{Proposition A3.} \emph{$T$ is a Plott function.}\medskip

{\em Proof}. Let us verify that $T$ satisfies the heritage and outcast properties.

For the heredity: let $A\subseteq B$ and  $a\in T (B)\cap A$.
Then, suppose on the contrary that $a\notin T (A)$, that is, $a\in L (A)$.
 If $a\in D$, then $a\in L (B)$, that      contradicts to the fact that $a\in T (B)$. If $a\notin D$, then $a\prec A$ and from L3 %RS
 $a\prec B$, $a\in L (B)$, a contradiction to $a\in T (B)$.

%\emph{Outcast.}
For the outcast: let $T (A)\subseteq B\subseteq A$. Then for a negligible set $A$, we have that $B$ is negligible as well and % we need to show that
 $T(A)=T(B)$. %The statement is obvious if $A$ is negligible, because then $B$  is negligible as well.
For an essential set $A$, from %We will assume that $A$ is essential. Since
$T(A)\subseteq B$ follows       $A\subseteq L(A)\cup A\subseteq L(A)\cup B$. Then, by  Lemma \ref{lemma2},  $L(A)\subseteq       L(B)$.
By the same Lemma applied to the inclusion $B\subseteq A$, we have $L(B)\subseteq L(A)$. From that we get  $L (A)=L (B)$ and $T (A)=T (B)$. $\Box$\medskip

Thus, the \emph{operations $G \mapsto \prec _G$ and $\prec \mapsto T$
are mutual inverse and establish a  bijection between the set of Plott functions and the set of Lehmann hyper-orders.}

\section*{Acknowledgments} We thank Arkady Slinko, the editor and referees for useful remarks and suggestions.

\end{document}